\documentclass[12pt,reqno]{amsart}

\usepackage{amsfonts,color,amsthm,amsmath,amssymb}

\def\Box{\vcenter{\vbox{\hrule\hbox{\vrule
     \vbox to 8.8pt{\hbox to 10pt{}\vfill}\vrule}\hrule}}}

\newcommand{\tr}{\textup{Tr}}
\newcommand{\Norm}{\textup{Norm}}
\def\qed{{\hfill$\square$}}
\def\mod{{\mathrm{mod\,\,}}}

\newcommand{\ga}{{\gamma}}
\newcommand{\la}{\langle}
\newcommand{\ra}{\rangle}

\newcommand{\F}{{\mathbb F}}
\newcommand{\Z}{{\mathbb Z}}

\newcommand{\Q}{{\mathbb Q}}

\newtheorem{thm}{Theorem}
\newtheorem{prob}[thm]{Problem}
\newtheorem{lem}[thm]{Lemma}

\newtheorem{conj}[thm]{Conjecture}
\newtheorem{remark}[thm]{Remark}
\numberwithin{equation}{section}

\begin{document}
\title{Nonsymmetric primitive translation schemes on prime power number of vertices}
\author{Tao Feng and Koji Momihara  }
\thanks{T. Feng is with Department of Mathematics, Zhejiang University, Hangzhou 310027, Zhejiang, China
(e-mail: tfeng@zju.edu.cn). T. Feng was supported  by
Fundamental Research Fund for the Central Universities of China, Zhejiang Provincial Natural Science Foundation under Grant LQ12A01019,  the National Natural Science Foundation of China under Grant 11201418, and  the Research Fund for Doctoral Programs from the Ministry of Education of China under Grant 20120101120089. }
\thanks{K. Momihara is with Faculty of Education, Kumamoto University,
2-40-1 Kurokami, Kumamoto 860-8555, Japan (e-mail:
momihara@educ.kumamoto-u.ac.jp). The work of
K. Momihara was supported by JSPS under Grant-in-Aid for Young Scientists (B) 25800093.}
\maketitle

\begin{abstract}
It is well-known that translation schemes on prime number of vertices are exactly the cyclotomic schemes. In this current paper, we show that there are no nonsymmetric primitive translation schemes on prime square vertices with at most four classes. On the other hand, we find new  non-symmetric four- and five-class association schemes from cyclotomy as fission schemes of certain symmetric three-class schemes. Moreover, we provide an affirmative answer to the following question raised by Song \cite{song_2}:
Are there any other two-class primitive schemes that admit symmetrizable fission schemes besides the cyclotomic scheme of index $2$  for $q \equiv5 \pmod{8}$?
To be more specific, we   show that a certain two-class primitive scheme  in the finite field $\F_{37^3}$ constructed by Feng and Xiang in \cite{fx} admits a four-class fission scheme. This fission scheme is realized as a fusion scheme of the cyclotomic scheme of index $28$.
\end{abstract}

\section{Introduction}

On p. 68 of the book \cite{bcn}, the following conjecture is suggested:
\begin{conj}\label{conj_bcn}
Let $(X,{\mathcal R})$ be a primitive translation association scheme. Then one of the following holds: (1) $X$ is elementary abelian, (2) $(X,{\mathcal R})$ is a Hamming scheme, (3) $(X,{\mathcal R})$ is of Latin square type or of negative Latin square type.
\end{conj}
A scheme is of (negative) Latin square type if each of its nontrivial relation is a strongly regular graph  of (negative) Latin square type. By A.V. Ivanov's classification, an amorphic association scheme with at least
three classes must be of Latin square type or negative Latin square type, cf. \cite{avi,amorphic}. He also conjectured in  \cite{aai} that if each nontrivial relation in an association scheme is strongly regular, then the association scheme must be amorphic. This conjecture turned out to be false, and all the known primitive counterexamples are constructed in finite fields using the union of cyclotomic classes; please refer to \cite{fwx} and the references therein for more details. In view of Conjecture \ref{conj_bcn} and these facts, finite fields and cyclotomy play important roles in the study of primitive association schemes by yielding new schemes which help us to understand their structures better.

By a well-known theorem due to Hanaki and Uno \cite{hu}, all association schemes with a prime number of points must be pseudocyclic, which is the analogue of the fact that finite groups of prime order must be cyclic. In the translation scheme case, we have the stronger result that a translation scheme with a prime number of points must be a cyclotomic scheme, which is a corollary of the following multiplier theorem of Schur and Wielandt.
\begin{thm}\cite[Thm. 2.10.1]{bcn}  Let $(X,\mathcal{R})$ be a translation scheme, and let $s$ be an integer coprime to $|X|$. Then for any $R\in\mathcal{R}$, we have $R^{(s)}=\{sx|\,x\in R\}\in\mathcal{R}$.
\end{thm}

In \cite[p. 58]{bi}, the question of when a  symmetric association scheme can be split into  non-symmetric commutative association scheme is posed. In \cite{ck}, the authors considered non-symmetric commutative association schemes with exactly one pair of nonsymmetric relations, and introduced a set of feasibility and
realizability conditions.
From  each nonsymmetric association scheme with three classes, a symmetric association scheme with two classes can be obtained by merging the nonsymmetric relations, which is called the symmetrization of the original scheme. Two-class skew-symmetric association schemes are regular tournaments, whose symmetrization is trivial. There are some papers on the feasibility conditions for the existence of such schemes, c.f.  \cite{bs,gc,s3}.  In \cite{lj}, the author did a systematic investigation of such schemes  with no more than $100$ vertices. There are several such schemes over $36$ vertices, see \cite{36_1,36_3,36_2,lj}. There is an infinite family of such schemes on $4s^4$ number of vertices where $s$ is a power of $2$, cf. \cite{3c_1,3c_2,ma_4}. In \cite{lj10}, the author found some more such schemes on $64$ vertices. These are all the known primitive three-class skew-symmetric schemes so far. On the nonexistence side, in \cite{gc} the authors showed that there is no three-class nonsymmetric scheme whose symmetrization is a Paley type strongly regular graph.

A skew-symmetric association scheme is a scheme with no symmetric adjacency relations other than the diagonal one. A typical example is given by the cyclotomic scheme of index $4$ in the finite field $\F_q$, where $q\equiv 5\pmod{8}$. The study of four-class skew symmetric schemes was initiated in \cite{song_2}. In \cite{bs}, Bannai and Song raised a question regarding the existence of four-class amorphous skew-symmetric association schemes, which was answered in the negative in \cite{ma}.   In the same paper, Ma classified four-class skew-symmetric schemes by their character tables, which fall into three types. In a subsequent paper \cite{mw}, the authors determined their intersection matrices, and generated a list of feasible parameters with small number of vertices. The aforementioned  cyclotomic scheme of index $4$ are the only known such schemes which are primitive to our knowledge.

In the current paper, we are primarily concerned with primitive translation association schemes on prime power number of points. Let $p$ be a  prime. Symmetric two-class translation schemes on $p^2$ vertices are exactly the strongly regular graphs of Latin square type with $\Z_p^*$ as multiplier group; see the survey \cite{pds_sur} for details. Nonsymmetric two-class translation schemes correspond to skew Hadamard difference sets whose ambient groups must have order congruent to $3$ modulo $4$, so there are no such schemes on $p^2$ vertices. We show that there are no three- or four-class primitive translation schemes on prime square vertices with at least one pair of nonsymmetric relations.
On the other hand, we find new  nonsymmetric four- and five-class association schemes   from cyclotomy as fissions of certain symmetric three-class schemes, which are primitive under certain conditions.
In particular, the obtained nonsymmetric four-class association schemes
have exactly one pair of nonsymmetric relations, which give
a lot of examples of association schemes studied in \cite{ck} by Chia and Kok.
Moreover, we provide an affirmative answer to the following question raised by Song \cite{song_2}:

\begin{prob}Are there any other two-class primitive schemes that admit symmetrizable fission schemes besides the cyclotomic scheme of index $2$  for $q \equiv5 \pmod{8}$?
\end{prob}
To be more specific, we shall show that a certain two-class primitive schemes in the finite field $\F_{37^3}$ constructed in \cite{fx} admit a four-class fission scheme. This fission scheme is realized as the fusion scheme of the cyclotomic scheme of index $28$.
\section{Preliminaries}\label{sec:pre}
Let $(X,\{R_i\}_{i=0}^d)$ be  a commutative $d$-class association scheme, and $A_i$ be the adjacency matrix of the relation $R_i$. Its Bose-Mesner algebra is defined as ${\mathcal A}=\langle A_0,A_1,\ldots,A_d\rangle$.    Denote by $E_0=(1/|X|)J,E_1,\ldots,E_d$ the primitive idempotents of ${\mathcal A}$, where $J$ is the $|X|\times |X|$ all one matrix. Both the $A_i$'s and the $E_j$'s form basis of the algebra ${\mathcal A}$, and we thus have
\[
A_j=\sum_{i=0}^{d}p_j(i)E_i,\, \, E_j=\frac{1}{|X|}\sum_{i=0}^dq_j(i)A_i,
\]
for some constants $p_j(i)$ and  $q_j(i)$. Let $P$ and $Q$ be the $(d+1)\times (d+1)$ matrix, whose $(i,j)$-th entry is $p_j(i)$ and $q_j(i)$ respectively. The matrix $P$ and $Q$ is called the {\it first eigenmatrix} (or the {\it character table}) and the
{\it second eigenmatrix} of  $(X,\{R_i\}_{i=0}^d)$ respectively.
There exists numbers $p_{i,j}^k$ and $q_{i,j}^k$ such that
\[
A_iA_j=\sum_{k=0}^d p_{i,j}^kA_k,\, \,  E_i\circ E_j=\sum_{k=0}^d q_{i,j}^kE_k,
\]
where $\circ$ stands for entrywise multiplication.
The numbers $p_{i,j}^k$'s and $q_{i,j}^k$'s are called the {\it intersection numbers} and
the {\it Krein parameters} of $(X,\{R_i\}_{i=0}^d)$, respectively. For each $0\leq i\leq d$, the $i$-th intersection matrix $B_i$ is defined as the $(d+1)\times (d+1)$ matrix whose $(k,j)$-th entry is $p_{i,j}^k$. Clearly, $B_0$ is the identity matrix.

Given two schemes ${\mathcal X}=(X,\{R_i\}_{i=0}^d)$ and
${\mathcal X}'=(X,\{R_i'\}_{i=0}^e)$, if for each $0\le i\le d$
$R_i\subseteq R_j'$ for some $0\le j\le e$, then we say that ${\mathcal X}$ is a {\it fusion scheme} of ${\mathcal X}'$, and ${\mathcal X}'$ is a {\it fission scheme} of ${\mathcal X}$.
We shall need the following well-known criterion due to Bannai \cite{BannaiSub} and Muzychuk \cite{Muzthesis}, called the {\it Bannai-Muzychuk criterion}: {\it Let $P$ be the first eigenmatrix of an association scheme $(X, \{R_i\}_{0\leq i\leq d})$, and $\Lambda_0:=\{0\}, \Lambda_1,\ldots ,\Lambda_{d'}$ be a partition of $\{0,1,\ldots ,d\}$. Then $(X, \{R_{\Lambda_i}\}_{0\leq i\leq d'})$ forms an association scheme if and only if there exists a partition $\{\Delta_i\}_{0\leq i\leq d'}$ of $\{0,1,\ldots ,d\}$ with $\Delta_0=\{0\}$ such that each $(\Delta_i, \Lambda_j)$-block of $P$ has a constant row sum. Moreover, the constant row sum of the $(\Delta_i, \Lambda_j)$-block is the $(i,j)$-th entry of the first eigenmatrix of the fusion scheme.}

We call an association scheme $(X,\{R_i\}_{i=0}^d)$ a {\it translation association scheme} or a {\it Schur ring} if $X$ is an  (additively written) finite abelian  group and there exists a partition $R_0'=\{0\},R'_1,\cdots,R'_d$  of $X$ such that
\[
R_i=\{(x,x+y)|\,x\in X, y \in R'_i\}.
\]
For brevity, we will just say that $(X,\{R'_i\}_{i=0}^d)$ is a translation scheme, and the $R'_i$'s are the relations. The set $R_0'$ defines the trivial relation and all the other relations are called nontrivial.

Assume that $(X,\{R_i\}_{i=0}^d)$ is a translation   scheme. There is an equivalence relation defined on the character group $\widehat{X}$ of $X$ as follows: $\chi\sim\chi'$ if and only if $\chi(R_i)=\chi'(R_i)$ for each $0\leq i\leq d$. Here $\chi(R)=\sum_{g\in R}\chi(g)$, for any $\chi\in\widehat{X}$, and $R\subseteq X$. Denote by $D_0, D_1,\cdots,D_d$ the equivalence classes, with $D_0$ consisting of only the principal character $\chi_0$. Then $(\widehat{X},\{D_i\}_{i=0}^d)$ forms a translation   scheme, called the {\it dual} of $(X,\{R_i\}_{i=0}^d)$. If $\{D_i\}_{i=0}^d$ is mapped to $\{R_i\}_{i=0}^d$ under the natural isomorphism between $X$ and $\widehat{X}$, then the scheme is called {\it self-dual}. The first eigenmatrix of the dual scheme is equal to the second eigenmatrix of the original scheme. Please refer to \cite{bi} and \cite{bcn} for more details.

An association scheme $(X,\{R\}_{i=0}^{d})$ is called {\it primitive} if each of the undirected graphs defined by $R_i$, $1\le  i\le d$, is connected, and {\it imprimitive} otherwise. A nonsymmetric scheme is primitive if and only if its symmetrization is primitive. If $(X,\{R_i\}_{i=0}^{d})$ is a symmetric translation scheme, then it is well-known that $(X,\{R_i\}_{i=0}^{d})$ is primitive if and only if 
each nontrivial relation $R_i$ generates the whole group $X$, i.e., not contained in any proper
subgroup of $X$; see \cite{bi}. It is equivalent to the fact that for each nontrivial relation $R_i$ ,
$\chi(R_i)\not=|R_i|$ for all nonprincipal characters $\chi$ of $X$; see \cite{bh}.

A classical example of translation schemes is the cyclotomic scheme which we describe now. Let $p$ be a prime and $q=p^f \,(f\geq 1)$ be a prime power, $N|q-1$, and $\gamma$ be a primitive root of the finite field $\F_q$. Define the multiplicative subgroup $C_0^{(N,q)}=\langle\gamma^N\rangle$. Its cosets $C_i^{(N,q)}=\gamma^iC_0^{(N,q)}$, $0\leq i\leq N-1$, are called the {\it cyclotomic classes} of order $N$ of $\F_q$.
Together with $\{0\}$, they form an $N$-class association scheme, which is called the {\it cyclotomic scheme} of index $N$. If $D=\bigcup_{i\in I}C_i^{(N,q)}$ for some subset $I\subseteq\Z_N$, then we say that $D$ is defined by the index set $I$.

Now we introduce some notations that we shall use in the later sections. We define the {\it Gauss periods} of index $N$ over $\F_q$ as follows:
\[
\eta_i=\sum_{x\in C_i^{(N,q)}}\psi(x),\;0\leq i\leq N-1,
\]
where $\psi$ is the canonical additive character of $\F_q$ defined by $\psi(x)=e^{\frac{2\pi i}{p}\tr_{q/p}(x)}$, $x\in \F_q$.
For each multiplicative character $\chi$ of $\F_q^\ast$, the multiplicative group of $\F_q$, we define the {\it Gauss sum}
\[
G_q(\chi)=\sum_{x\in \F_q^\ast}\psi(x)\chi(x).
\]
The following relation will be repeatedly used in this paper (cf. \cite[P.~195]{LN97}):
\[
\psi(x)=\frac{1}{q-1}\sum_{\chi\in\widehat{\F_q^\ast}}G_q(\chi)\chi^{-1}(x),\;\forall\, x\in \F^\ast,
\]
where $\widehat{\F_q^\ast}$ is the character group of $\F_q^\ast$.
The Gauss period can be expressed as a linear combination of Gauss sums as follows:
\begin{align*}
\eta_i&=\psi(\omega^i C_0^{(N,q)})
=\frac{1}{q-1}\sum_{\chi\in\widehat{\F_q^\ast}}G_{q}(\chi)\chi^{-1}(\ga^i)\sum_{x\in C_0^{(N,q)}}\chi^{-1}(x)\\
&=\frac{1}{N}
\sum_{i=0}^{N-1}G_q(\varphi^{-i})\varphi(\ga^{i} ),
\end{align*}
where $\varphi$ is a multiplicative character of order $N$ of $\F_q^\ast$.

We record some basic properties of Gauss sums here \cite{LN97}:
\begin{itemize}
\item[(i)] $G_q(\chi)\overline{G_q(\chi)}=q$ if $\chi$ is nontrivial;
\item[(ii)] $G_q(\chi^p)=G_q(\chi)$, where $p$ is the characteristic of $\F_q$;
\item[(iii)] $G_q(\chi^{-1})=\chi(-1)\overline{G_q(\chi)}$;
\item[(iv)] $G_q(\chi)=-1$ if $\chi$ is trivial.
\end{itemize}
In general, explicit evaluations of Gauss sums are very difficult. There are only a few cases where
the Gauss sums have been evaluated, see \cite{BEW97} for known results on Gauss sums.
In Section~\ref{secind2}, we will use the explicit evaluations of Gauss sums of  the quadratic residue case and index 2 case.
\begin{thm}(\cite{LN97})\label{le:quad}
Let $\eta$ be the quadratic character of $\F_q=\F_{p^f}$, and $p^\ast=(-1)^{\frac{p-1}{2}}p$. It holds that
\[
G_q(\eta)=(-1)^{f-1}\left(\sqrt{p^*}\right)^f.
\]
\end{thm}
\begin{thm}\label{Sec2Thm2}(\cite[Case D; Theorem~4.12]{YX10})
Let $N=2p_1^m$, where $p_1>3$ is a prime such that $p_1\equiv 3\,(\mod{4})$ and $m$ is a positive integer. Assume that $p$ is a
prime such that  $[\Z_N^\ast:\langle p \rangle]=2$.
Let $f=\phi(N)/2$, $q=p^f$, and $\chi$ be a multiplicative character of order $N$ of $\F_q$. For
$0\le t\le m-1$, we have
\begin{eqnarray*}
G_q(\chi^{p_1^t})&=&\left\{
\begin{array}{ll}
(-1)^{\frac{p-1}{2}(m-1)}p^{\frac{f-1}{2}-hp_1^t}
\sqrt{p^\ast}\left(\frac{b+c\sqrt{-p_1}}{2}\right)^{2p_1^t}&  \mbox{if $p_1\equiv 3\,(\mod{8})$,}\\
(-1)^{\frac{p-1}{2}m}p^{\frac{f-1}{2}}
\sqrt{p^\ast}&  \mbox{if $p_1\equiv 7\,(\mod{8})$;}
 \end{array}
\right.\\
G_q(\chi^{2p_1^t})&=&p^{\frac{f-p_1^t h}{2}}\left(\frac{b+c\sqrt{-p_1}}{2}\right)^{p_1^t};\\
G_q(\chi^{p_1^m})&=&(-1)^{\frac{p-1}{2}\frac{f-1}{2}}p^{\frac{f-1}{2}}\sqrt{p^\ast},
\end{eqnarray*}
where  $h$ is the class number of $\Q(\sqrt{-p_1})$, and $b$ and $c$ are integers
determined by $4p^{h}=b^2+p_1c^2$ and $bp^{\frac{f-h}{2}}\equiv -2\,(\mod{p_1})$.
\end{thm}
We shall also need the following  {\it Davenport-Hasse lifting formula} on Gauss sums.
\begin{thm}\label{thm:lift}(\cite[Theorem~5.14]{LN97})
Let $\chi$ be a nontrivial multiplicative character of $\F_q=\F_{p^f}$ and
let $\chi'$ be the lifted character of $\chi$ to the extension field $\F_{q'}=\F_{p^{fs}}$, i.e., $\chi'(\alpha):=\chi(\Norm_{\F_{q'}/\F_q}(\alpha))$ for any $\alpha\in \F_{q'}^\ast$.
It holds that
\[
G_{q'}(\chi')=(-1)^{s-1}(G_q(\chi))^s.
\]
\end{thm}

\section{Nonsymmetric primitive schemes on prime square number of vertices}
In this section, we show that there is no  nonsymmetric primitive translation scheme with at most four classes on prime square number of vertices.  As we mentioned in the introduction, there are no such schemes with two classes. By the Wielandt-Schur multiplier theorem, such schemes can only survive in the additive group $G=(\F_{p^2},+)$. Let $(G,\{R_i\}_{i=0}^d)$ be such a scheme with $d=3, 4$. By examining the maximal symmetric subschemes and using the fact that the dual of a symmetric scheme is still symmetric, we see that $(G,\{R_i\}_{i=0}^d)$ has the same number of nonsymmetric relations as its dual. Write $\ga$ for a fixed primitive root of $\F_{p^2}$, and $\ga_0=\ga^{p+1}$ in $\F_p^\ast$.

\begin{lem}\label{mult}
The scheme $(G,\{R_i\}_{i=0}^d)$ has the set $S$ of nonzero squares of $\Z_p^\ast$ as a multiplier group.
\end{lem}
\proof
This follows from Wielandt-Schur multiplier theorem. Let $R$ be any nontrivial relation. If $R^{(\ga_0^2)}\ne R$, then $R,R^{(\ga_0)},R^{(\ga_0^2)},R^{(\ga_0^3)}$ are the four distinct nontrivial relations of the scheme. It follows that the scheme has two pairs of nonsymmetric relations, and $R^{(-1)}=R^{(\ga_0^2)}$. Its symmetrization is a two-class scheme $(G, \{\{0\}\}\cup \{R\cup R^{(\ga_0^2)},R^{(\ga_0)}\cup R^{(\ga_0^3)}\})$, so $R\cup R^{(\ga_0^2)}$ is a partial difference set in $G$ which has $\Z_p^\ast$ as a multiplier group, cf. \cite{pds_sur}. It contradicts to the fact that $R\cup R^{(\ga_0^2)}\ne R^{(\ga_0)}\cup R^{(\ga_0^3)}$.
\qed
\vspace{0.3cm}

One consequence of the above lemma is that $-1$ is a nonsquare in $\Z_p$, and so $p\equiv 3\pmod{4}$. It is clear that $S=C_0^{(N,p^2)}$ with $N=2(p+1)$, the multiplicative subgroup of index $N=2(p+1)$ in $\F_{p^2}$. Let $C_i:=\ga^iS$, $0\leq i\leq N-1$ be the cyclotomic classes of index $N$. If a set $D=\bigcup_{i\in I}C_i $, then we say that $D$ is defined by the index set $I\subset\Z_N$.  According to Lemma \ref{mult}, we assume that the nontrivial relations in the scheme $(G,\{R_i\}_{i=0}^d)$ are defined by the index sets $I_1,\ldots,I_d$, and the nontrivial relations in the dual scheme are defined by the index sets $J_1,\ldots,J_d$. Since $-1\in C_{p+1}$, we have $J_i=\{i+p+1\pmod{N}\,|\,i\in J_i\}$, and $J_i([0]-[p+1])=0$ in $\Z[\Z_N]$ if $J_i$ defines a symmetric relation.

Throughout this section, we use the standard notations on group rings as in \cite{pott}. To make distinction, we use $[i]$ for the group ring element in  $\Z[\Z_{N}]$ corresponding to $i\in\Z_{N}$. Let $L=\{x\in\F_{p^2}^*\,|\,\tr(x)=1\}$ which has size $p$, where $\tr:=\tr_{p^2/p}$. We observe that $T_s=\{i\pmod{N}\,|\,\gamma^i\in L\}\subset \Z_N$ is a relative difference set in $\Z_N$ w.r.t. the order $2$ subgroup $\la N_0\ra$ (cf. \cite{pott}):
\[
T_sT_s^{(-1)}=p+\frac{p-1}{2}(\Z_N-\la N_0\ra).
\]
For a group ring element $Y\in \Z[\Z_N]$ and $i\in\Z_N$, we use $[Y]_i$ to stand for the coefficient of $i$ in $Y$.

Let $\psi$ be the canonical additive character of $\F_{p^2}$ such that  $\psi(x)=\xi_p^{\tr(x)}$, where $\xi_p$ is a primitive $p$-th root of unity. By direct computations, we have
\begin{align*}
\psi(C_i) =\sum_{j=0}^{M-1}\xi_p^{\tr(\gamma^{i+Nj})}
=\sum_{j=0}^{M-1}\xi_p^{ (\ga^{Nj}\tr (\gamma^{i}))}
=\psi(\tr (\gamma^{i})S ). 
\end{align*}

We have a partition of $\Z_N$ into three parts by $T_0=\{\frac{p+1}{2},\frac{3(p+1)}{2}\}$, $T_s$ and $T_n:=\{i+(q+1)\pmod{N}\,|\,i\in T_s\}$. If $i\in  T_0$, then $\tr (\ga^i)=0$, and $\psi(C_i)=|C_0|$.   If $i\in T_s$, then  $\tr (\ga^i)$ is a nonzero square, and $\psi(C_i)=\frac{-1+\sqrt{p}}{2}$ with a proper choice of $\xi_p$. If $i\in T_n$, then  $\tr (\ga^i)$ is a nonsquare, and $\psi(C_i)=\frac{-1-\sqrt{p}}{2}$. Assume that  the set $D$   is defined by the index set $I\subset \Z_N$. Then we have
\begin{align*}
\psi(\ga^aD)&=\sum_{i\in I}\psi(C_{i+a})
=M\cdot|I+a\cap T_0|+\frac{-1+\sqrt{p}}{2}|I+a\cap T_s|+\frac{-1-\sqrt{p}}{2}|I+a\cap T_n|\\
&=\left[\left(\frac{p-1}{2}T_0-\frac{1}{2}T_s-\frac{1}{2}T_n+\frac{\sqrt{p}}{2}(T_s-T_n)\right)I^{(-1)}\right]_a\\
&=\frac{\sqrt{p}}{2}\left[(\sqrt{p}\cdot T_0+T_s-T_n)I^{(-1)}\right]_a-\frac{|I|}{2}.
\end{align*}

In the case where $D$ is a relation in the scheme $(G,\{R_i\}_{i=0}^d)$, the group ring element $(\sqrt{q}\cdot T_0+T_s-T_n)I^{(-1)}=\sum_{i=1}^dc_iJ_i$ for some algebraic integers $c_i$, $1\leq i\leq d$. Multiplying both sides with $[0]-[p+1]$, we have
$T_sI^{(-1)}([0]-[p+1])=\sum_{i}c_iJ_i([0]-[p+1])$, where the summation is over such $i$ that $J_i$ defines a nonsymmetric relation.

Now assume that the scheme $(G,\{R_i\}_{i=0}^d)$ has exactly one pair of nonsymmetric relations, defined by $I_1$ and $I_2$ respectively. The above argument yields that $T_sI_1^{(-1)}([0]-[p+1])= cJ([0]-[p+1])$ where the  $J$ defines a nonsymmetric relation in the dual scheme. Applying the same argument to the dual scheme, we get $T_s J^{(-1)}([0]-[p+1])= dI_1([0]-[p+1])$ for some $d\in\Z$. Multiplying both sides by $T_s^{(-1)}$ and taking involution, we get $T_sI_1^{(-1)}([0]-[p+1])=\frac{p}{d} J([0]-[p+1])$. Therefore, $p=cd$.  W.l.o.g., we assume that $c=\pm p$. Now $|T_s|=p$, so  $T_sI_1^{(-1)}$ has all coefficients  not exceeding $p$, and the coefficients $\pm p$ occur in $T_s I_1^{(-1)}([0]-[p+1])$ only if $|I_1|=p$ in which case our scheme is clearly imprimitive. Therefore, this case can not occur.

Now assume that $d=4$ and $(G,\{R_i\}_{i=0}^4)$ has two pair of nonsymmetric relations, defined by $I,I+(p+1),J,J+(p+1)\in\Z_N$ respectively. Also assume that the nontrivial relations of the dual scheme are defined by $I',I'+(p+1),J',J'+(p+1)\in\Z_N$. Now we apply the same argument as above and get
\begin{align*}
T_sI^{(-1)}([0]-[p+1])=aI'([0]-[p+1])+bJ'([0]-[p+1]),\\
T_sJ^{(-1)}([0]-[p+1])=cI'([0]-[p+1])+dJ'([0]-[p+1]).
\end{align*}
Multiplying both sides with $T_s^{(-1)}$ and taking involution, we have
\begin{align*}
&pI([0]-[p+1])=aT_sI'^{(-1)}([0]-[p+1])+bT_sJ'^{(-1)}([0]-[p+1]),\\
&pJ([0]-[p+1])=cT_sI'^{(-1)}([0]-[p+1])+dT_sJ'^{(-1)}([0]-[p+1]).
\end{align*}
The determinant $ad-bc\ne 0$, since the left hand side does not differ by a constant multiple. We may solve that
\begin{align*}
&T_sI'^{(-1)}([0]-[p+1])=\frac{p}{ad-bc}(dI-bJ)([0]-[p+1]),\\
&T_sJ'^{(-1)}([0]-[p+1])=\frac{p}{ad-bc}(-cI+aJ)([0]-[p+1]).
\end{align*}
Therefore, $ad-bc$ divides $p\cdot\gcd(b,d,a,c)$. Since $\gcd(b,d,a,c)^2|ad-bc$, we have $\gcd(b,d,a,c)|p$. It follows that $\gcd(b,d,a,c)=1$ or $p$. It is easy to see that each of these four numbers are less than $p$, so we have $\gcd(b,d,a,c)=1$, and $ad-bc\in\{\pm1,\pm p\}$. By replacing the scheme with its dual if necessary, we may assume that $ad-bc\in\{\pm 1\}$.

Now we solve that
\begin{align*}
 &I'^{(-1)}([0]-[p+1])=\frac{1}{ad-bc}T_s^{(-1)}(dI-bJ)([0]-[p+1]),\\
 &J'^{(-1)}([0]-[p+1])=\frac{1}{ad-bc}T_s^{(-1)}(-cI+aJ)([0]-[p+1]).
\end{align*}
Therefore, multiplying each equation by its involution respectively, we get
\begin{align*}
 &I'I'^{(-1)}([0]-[p+1])= p(dI-bJ)(dI-bJ)^{(-1)}([0]-[p+1]),\\
 &J'J'^{(-1)}([0]-[p+1])= p(-cI+aJ)(-cI+aJ)^{(-1)}([0]-[p+1]).
\end{align*}
Assume w.l.o.g. that $1<|I'|\leq\frac{p+1}{2}$ ($|I'|=1$ would imply imprimitivity). Observe that $ I'I'^{(-1)}$ have coefficients between $0$ and $|I'|$, so  $I'I'^{(-1)}([0]-[p+1])$ has coefficient between $- |I'|$ and $|I'|$. Therefore, together with the fact that it is divisible by $p$, we have  $I'I'^{(-1)}([0]-[p+1])=0$. It follows that the coefficient of $p+1$ in  $I'I'^{(-1)}$ is the same as that of $0$, namely $|I'|$. This means $I'=I'+(p+1)$, which  is a contradiction. This completes the proof.\\


We summarize the results in this section in the following theorem.
\begin{thm}There is no nonsymmetric primitive translation scheme over prime square number of vertices with at most four classes.
\end{thm}
\section{Fission schemes of three-class association schemes based on
index $2$ Gauss sums }\label{secind2}
Throughout this section, we assume that we are in the same settings
as in  Theorem~\ref{Sec2Thm2}: $p_1>3$ is a prime such that $p_1\equiv 3\pmod{4}$; $p$ is an odd prime such that  $[\Z_{2p_1}^\ast:\langle p\rangle]=2$ (that is, $f:={\rm ord}_{2p_1}(p)=(p_1-1)/2$). Put $q=p^f$.
In this case,
the Gauss sums $G_{q}(\chi^j)$ were completely evaluated in \cite{YX10} as described in Theorem~\ref{Sec2Thm2},
where $\chi$ is a multiplicative character of order $2p_{1}$ of $\F_q$.
(We called this the {\it index $2$ Gauss sum}.)

Consider the extension field $\F_{q^s}$ for a positive integer $s$ and put
$G=(\F_{q^s},+)$. Let $\gamma$ be a primitive root of $\F_{q^s}$. Consider the automorphism group of $G$ generated by multiplication by $\ga^{p_1}$ and the Frobenius automorphism which raises each element of $\F_{q^s}$ to its $p$-th power. Its orbits on $(\F_{q^s},+)$ consist of
\[
R_0=\{0\},\, \, R_1=\bigcup_{j \in \langle p\rangle (\mod{p_1})}C_{j}^{(p_1,q^s)}, \,
\, R_2=\bigcup_{j \in -\langle p\rangle (\mod{p_1})}C_{j}^{(p_1,q^s)},\,
R_3=C_{0}^{(p_1,q^s)}.
\]
It is clear that this partition $(G,\{R_i\}_{i=0}^{3})$ becomes
a  symmetric three-class self-dual scheme. This scheme is primitive if and only if $R_3$ generates the whole group, namely $\F_p[\ga^{p_1}]=\F_{q^s}$. This happens exactly when $\ga^{p_1}$ has a minimal polynomial of degree $fs$ over $\F_p$, i.e., $p$ has order $fs$ modulo $\frac{q^s-1}{p_1}$.

We now consider a four-class fission scheme of this scheme.
\begin{thm}\label{ind2} Take the same notations as above, and further assume that $p_1\equiv 7\,(\mod{8})$. Define
\[
S_0=R_0, \, S_1=R_1, \, S_2=R_2,\,
S_3=C_{0}^{(2p_1,q^s)},\,
S_4=C_{p_1}^{(2p_1,q^s)}. 
\]
Then $(G,\{S_i\}_{i=0}^{4})$ becomes
a four-class self-dual association scheme.
\end{thm}
\proof
In view  of the Bannai-Muzychuk criterion, we need to compute
the eigenvalues $\psi(\gamma^a S_i)$, $1\le i\le 4$, and show that they  are constants according as $\gamma^a \in S_j$, $1\le j\le 4$. Since
we already know that  $\psi(\gamma^a S_i)$, $i=1,2$, are  constants according as $\gamma^a \in S_j$, $1\le j\le 4$, we
compute only the sum $\psi(\gamma^a S_3)$.
 (Note that $\psi(\gamma^a S_4)=\psi(\gamma^a R_3)-\psi(\gamma^a S_3)$.)  By the expression of Gauss periods as linear combinations of Gauss sums as explained in the introduction, we have
\begin{align}\label{S_3}
\psi(\gamma^a S_3)&=\frac{1}{2p_1}\sum_{i=0}^{2p_1-1}G_{q^s}(\chi^i)\chi^{-i}(\gamma^a)\nonumber\\
&=\frac{1}{2p_1}\sum_{0\le i\le 2p_1-1: \,i\,odd}^{}G_{q^s}(\chi^{i})\chi^{-i}(\gamma^a)+\frac{1}{2p_1}\sum_{0\le i\le 2p_1-1: \,i\,even}^{}G_{q^s}(\chi^{i})\chi^{-i}(\gamma^a),
\end{align}
where $\chi$ is a multiplicative character of order $2p_1$ of $\F_{q^s}$.
On the other hand, we have
\[
\psi(\gamma^a R_3)=\frac{1}{p_1}\sum_{i=0}^{p_1-1}G_{q^s}(\chi^{2i})\chi^{-2i}(\gamma^a),
\]
which implies that
\[
\eqref{S_3}=\frac{1}{2p_1}\sum_{0\le i\le 2p_1-1: \,i\,odd}^{}G_{q^s}(\chi^{i})\chi^{-i}(\gamma^a)+\frac{1}{2}\psi(\gamma^a R_3).
\]
By Theorem~\ref{Sec2Thm2} and the Davenport-Hasse lifting formula,
\[
G_{q^s}(\chi^i)=(-1)^{s-1}(-1)^{\frac{p-1}{2}s}p^{\frac{f-1}{2}s}
\sqrt{p^\ast}^s
\]
for all odd $i$, where we used the facts that  $\chi^i(-1)\overline{\sqrt{p^\ast}}=\sqrt{p^\ast}$ and $\frac{p-1}{2}\equiv \frac{p-1}{2}\frac{f-1}{2}\,(\mod{2})$.
We now compute that
\begin{align*}
\sum_{0\le i\le p_1-1: \,i\,odd}G_{q^s}(\chi^{i})\chi^{-i}(\gamma^a)&=
(-1)^{s-1}(-1)^{\frac{p-1}{2}s}p^{\frac{f-1}{2}s}
\sqrt{p^\ast}^s\chi^{-1}(\gamma^a)\sum_{i=0}^{p_1-1}\chi^{-2i}(\gamma^a)\\
&=
(-1)^{s-1}(-1)^{\frac{p-1}{2}s}p^{\frac{f-1}{2}s}
\sqrt{p^\ast}^s\cdot \left\{
\begin{array}{ll}
(-1)^{a}&  \mbox{if $a\equiv 0\,(\mod{p_1})$,}\\
0&  \mbox{otherwise.}
 \end{array}
\right.
\end{align*}
This implies that $\psi(\gamma^a S_3)$ is  constant according as $\gamma^a\in S_i$, $1\le i\le 4$, since we already know that $\psi(\gamma^a R_3)$ is a constant
according as $\gamma^a\in R_i$, $1\le i\le 3$. This completes the proof. 
\qed

\begin{remark}\label{primitivity1}
\begin{enumerate}
 \item If $s$ is odd and $p\equiv 3\,(\mod{4})$, then the resulting schemes of Theorem~\ref{ind2} are nonsymmetric four-class scheme with exactly one pair of nonsymmetric relations. This result gives a lot of  examples of nonsymmetric  association schemes studied in \cite{ck} by Chia and Kok.
 \item We investigate the primitivity of
     the resulting nonsymmetric association schemes in the case $s=1$ in  detail: As described in Section~\ref{sec:pre}, the nonsymmetric association scheme $(G,\{S_i\}_{i=0}^{4})$ is primitive if and only if  $\psi(\gamma^a(S_i\cup -S_i))\not=|S_i\cup -S_i|$ for $1\le i\le 4$.  Noting that $S_3\cup -S_3=R_3$, it is enough to see that $\psi(\gamma^a R_3)\not=|R_3|$. (Note that if the graph corresponding to $R_3$ is connected, then so are the others.) In \cite{FX111}, the sum $\psi(\gamma^a R_3)$ was computed as follows:
\[
\psi(\gamma^a R_3)=\frac{1}{p_1}\left(p^{(f-h)/2}\left(\frac{-b\pm p_1c}{2}\right)-1\right)\, \, \mbox{ or }\,\, \frac{1}{p_1}\left(p^{(f-h)/2}\left(\frac{p_1-1}{2}\right)-1\right).
\]
Therefore, $(G,\{S_i\}_{i=0}^{4})$ is primitive if and only if
$-b\pm p_1 c\not=2p^{\frac{f+h}{2}}$ and
$p_1-1\not=2p^{\frac{f+h}{2}}$.
Noting that $b^2+p_1c^2=4p^h$, for each $x\in \{-b\pm p_1c,p_1-1\}$ we have
$x\le b^2+c^2 p_1=4p^h< 2p^{(f+h)/2}$ if  $p_1>  2h+1$ and $c\not=0$.
Thus, the nonsymmetric association scheme $(G,\{S_i\}_{i=0}^{4})$ is primitive if $p_1> 2h+1$  and $c\not=0$.
\end{enumerate}
\end{remark}

Next we describe an infinite family of nonsymmetric  five-class association schemes with exactly two pairs of nonsymmetric relations in the case $p_1\equiv 3\,(\mod{8})$. Again, we  assume that we are in the same settings
as in  Theorem~\ref{Sec2Thm2} with $p_1\equiv 3\,(\mod{8})$. Let
$f=(p_1-1)/2$, $q=p^f$, and
$G=(\F_q,+)$. Define
\[
R_0=\{0\},\, \,
R_1=\bigcup_{j \in \langle p\rangle (\mod{p_1})}C_{j}^{(p_1,q)}, \,
\, R_2=\bigcup_{j \in -\langle p\rangle (\mod{p_1})}C_{j}^{(p_1,q)},\, \,
R_3=C_{0}^{(p_1,q)}.
\]
It is routine to show that $(G,\{R_i\}_{i=0}^{3})$ is
a three-class self-dual scheme.  This scheme is primitive if and only if $R_3$ generates the whole group, namely $\F_p[\ga^{p_1}]=\F_{q}$. This happens exactly when $\ga^{p_1}$ has a minimal polynomial of degree $f$ over $\F_p$, i.e., $p$ has order $f$ modulo $\frac{q-1}{p_1}$.

Suppose it holds  that $1+p_1=4p^h$, where $h$ is the class number of $\Q(\sqrt{-p_1})$.
Let $\xi_{q-1}=e^{\frac{2\pi i}{q-1}}$ and $\mathfrak{P}$ be a
prime ideal in $\Z[\xi_{q-1}]$ lying over $p$. Then,
$\Z[\xi_{q-1}]/\mathfrak{P}$ is the finite field of order $q$ and
\[
\Z[\xi_{q-1}]/\mathfrak{P}=\{\overline{\xi}_{q-1}^i\,|\,0\le i\le q-2\}\cup \{\overline{0}\},
\]
 where
$\overline{\xi}_{q-1}=\xi_{q-1}+\mathfrak{P}$.
Hence, $\gamma:=\overline{\xi}_{q-1}$ is a primitive root of
$\F_q=\Z[\xi_{q-1}]/\mathfrak{P}$.
Let $\omega_{\mathfrak{P}}$ be the
Teichm\"{u}ller character  of $\F_q^*$ such that
$\omega_{\mathfrak{P}}(\gamma)=\xi_{q-1}$.
Put $\chi:=\omega_{\mathfrak{P}}^\frac{q-1}{2p_1}$. Then
$\chi$ is a multiplicative character of order $2p_1$ of $\F_q^*$.

By Theorem~\ref{Sec2Thm2}, we have
\begin{equation}\label{Subsec2Eq1}
G_q(\chi)=p^{\frac{f-1}{2}-h}
\sqrt{p^\ast}\left(\frac{b+c\sqrt{-p_1}}{2}\right)^{2},
\end{equation}
where $b,c\not \equiv 0\,(\mod{p})$, $b^2+c^2p_1=4p^h$, and
$bp^\frac{f-h}{2}\equiv -2\,(\mod{p_1})$. Since
$1+p_1=4p^h$, we have $b,c\in \{-1,1\}$, where the sign of $c$ depends on the choice of $\mathfrak{P}$. It was shown in \cite{fx} that $bc\equiv -\sqrt{-p_1}\,(\mod{\mathfrak{P}})$.
On the other hand,   $1+p_1=4p^h$ implies that
$(1+\sqrt{-p_1})(1-\sqrt{-p_1})\in \mathfrak{P}$, so either $1+\sqrt{-p_1}\in \mathfrak{P}$ or
$1-\sqrt{-p_1}\in \mathfrak{P}$. We choose a prime ideal $\mathfrak{P}$ such that $1+\sqrt{-p_1} \in  \mathfrak{P}$.
Then  $bc\equiv -\sqrt{-p_1}\,(\mod{\mathfrak{P}})$  yields that $bc=1$.
From now on, we fix this choice of $\mathfrak{P}$ and the corresponding character $\chi$.

\begin{thm}\label{thmind22} Take the same notations as above.
Define
\begin{align*}
S_0=R_0, \, S_1=R_1,\, S_2=\bigcup_{j \in -\langle p\rangle (\mod{2p_1})}C_{j}^{(2p_1,q)}, \,S_3=\bigcup_{j \in -\langle p\rangle (\mod{2p_1})}C_{j+p_1}^{(2p_1,q)},
S_4=C_{0}^{(2p_1,q)},\,
S_5=C_{p_1}^{(2p_1,q)}.
\end{align*}
Then
$(G,\{S_i\}_{i=0}^5)$ becomes a five-class self-dual scheme.
\end{thm}
\proof Take the same notations as introduced above.
Note that $S_2\cup S_3=R_2$ and $S_4 \cup S_5=R_3$.
By the  Bannai-Muzychuk criterion, it is enough to show that the sums $\psi(\gamma^a S_2)$ and  $\psi(\gamma^a S_4)$
 are constants according as $\gamma^a\in S_i$, $1\le i\le 5$.

First, we compute  the sum $\psi(\gamma^a S_4)$.
By the same argument as that in the proof of Theorem \ref{ind2}, we have
\begin{align*}
\psi(\gamma^a S_4)&=\frac{1}{2p_1}\sum_{\ell=0}^{2p_1-1}G_{q}(\chi^\ell)\chi^{-\ell}(\gamma^{a} )=\frac{1}{2p_1}\sum_{0\le \ell\le 2p_1-1:\,\ell\, odd}G_{q}(\chi^{\ell})\chi^{-\ell}(\gamma^{a})+\frac{1}{2}\psi(\gamma^a R_3),
\end{align*}
where $\chi$ is as defined above.
Since $\psi(\gamma^a R_3)$ is constant
according as $\gamma^a\in R_i$, $1\le i\le 3$, we only need to consider the sum
\begin{equation}\label{eqn_2}
\sum_{0\le \ell\le 2p_1-1:\,\ell\, odd}G_{q}(\chi^{\ell})\chi^{-\ell}(\gamma^{a}).
\end{equation}
By Theorem~\ref{Sec2Thm2}, we have
\[
G_q(\chi^\ell)=A\left(\frac{b+c\sqrt{-p_1}}{2}\right)^2
\]
for all $\ell\in \langle p\rangle$, where $b,c$ are the same as in the
evaluation (\ref{Subsec2Eq1}) of $G_q(\chi)$ and
$A=p^{\frac{f-1}{2}-h}
\sqrt{p^\ast}$. By the choice of $\mathfrak{P}$,
it is expanded  as
\[
G_q(\chi^\ell)=A\left(\frac{1-p_1+2\sqrt{-p_1}}{4}\right).
\]
Since $G_q(\chi^{p_1})=p^{\frac{f-1}{2}}\sqrt{p^\ast}$ by Theorem~\ref{Sec2Thm2} and $\chi^\ell (-1)\overline{\sqrt{p^\ast}}=\sqrt{p^\ast}$ for any odd
$\ell$,
the sum in Eqn. \eqref{eqn_2} is reformulated as
\begin{equation}\label{eq11}
A\left(\left(\frac{1-p_1+2\sqrt{-p_1}}{4}\right)\sum_{\ell\in \langle p\rangle\,(\mod{2p_1})}
\chi^{-\ell}(\gamma^{a})+\left(\frac{1-p_1-2\sqrt{-p_1}}{4}\right)\sum_{\ell\in -\langle p\rangle\,(\mod{2p_1})}
\chi^{-\ell}(\gamma^{a})+p^h (-1)^a\right).
\end{equation}
Let $\eta$ be the quadratic character of $\F_{p_1}$ and
$\psi_{p_1}$ be the canonical additive character of $\F_{p_1}$.
Noting that $2$ is a nonsquare in $\F_{p_1}$,  it holds that
\begin{eqnarray*}
\sum_{\ell\in\langle p\rangle (\mod{2p_1})}\chi^{-\ell}(\gamma^a)&=&
\sum_{\ell\in\langle p\rangle (\mod{2p_1})}\chi^{-p_1 \ell}(\gamma^a)\chi^{2\frac{p_1-1}{2}\ell}(\gamma^a)\\
&=&
(-1)^a\frac{1}{2}\sum_{\ell\in \F_{p_1}^\ast}(1+\eta(\ell))\psi_{p_1}(-2^{-1}a\ell)\\
&=&
(-1)^a\frac{-1+\eta(-2^{-1}a)G_{p_1}(\eta)}{2}=(-1)^a
\frac{-1+\eta(a)\sqrt{-p_1}}{2}.
\end{eqnarray*}
Hence, by using $1+p_1=4p^h$, we have
\begin{align*}
\eqref{eq11}=&A(-1)^a\cdot \left\{
\begin{array}{ll}
0&  \mbox{if $a \in \langle p\rangle\,(\mod{p_1})$,}\\
\frac{1}{4}(-1+3p_1)+p^h&  \mbox{if $a \in -\langle p\rangle\,(\mod{p_1})$,}\\
-\frac{(p_1-1)^2}{4}+p^h&  \mbox{if $a \equiv 0\,(\mod{p_1})$.}
 \end{array}
\right.
\end{align*}
This shows that $\psi(\gamma^a S_4)$
 is constant according as $\gamma^a\in S_i$, $1\le i\le 5$.

Next, we consider the sum $\psi(\gamma^a S_2)$. In the same way, we get
\begin{align*}
\psi(\gamma^a S_2)&=\frac{1}{2p_1}\sum_{\ell=0}^{2p_1-1}\sum_{i\in -\langle p\rangle(\mod{2p_1})}G_{q}(\chi^\ell)\chi^{-\ell}(\gamma^{a+i} )\\
&=\frac{1}{2p_1}\sum_{0\le \ell\le 2p_1-1:\,\ell\, odd}\sum_{i\in -\langle p\rangle(\mod{2p_1})}G_{q}(\chi^{\ell})\chi^{-\ell}(\gamma^{a+i})+
\frac{1}{2}\psi(\gamma^a R_2).
\end{align*}
If $a\not\equiv 0\,(\mod{p_1})$, then
\begin{align*}
&\sum_{0\le \ell\le 2p_1-1:\,\ell\, odd}\sum_{i\in -\langle p\rangle (\mod{2p_1})}G_{q}(\chi^{\ell})\chi^{-\ell}(\gamma^{a+i})\\
=&A\left(\left(\frac{1-p_1+2\sqrt{-p_1}}{4}\right)\sum_{i\in -\langle p\rangle(\mod{2p_1})}\sum_{\ell\in \langle p\rangle(\mod{2p_1})}
\chi^{-\ell}(\gamma^{a+i})\right.\\
&\hspace{2cm}+\left.\left(\frac{1-p_1-2\sqrt{-p_1}}{4}\right)\sum_{i\in -\langle p\rangle(\mod{2p_1})}\sum_{\ell\in -\langle p\rangle (\mod{2p_1})}
\chi^{-\ell}(\gamma^{a+i})-\frac{p_1-1}{2}p^h (-1)^a\right)\\
=&A(-1)^a\left(\left(\frac{1-p_1+2\sqrt{-p_1}}{4}\right)\left(\frac{-1+\eta(a)\sqrt{-p_1}}{2}\right)\sum_{i\in \langle p\rangle (\mod{2p_1})}
\chi^{i}(\gamma)\right.\\
&\hspace{2cm}+\left.\left(\frac{1-p_1-2\sqrt{-p_1}}{4}\right)
\left(\frac{-1-\eta(a)\sqrt{-p_1}}{2}\right)\sum_{i\in -\langle p\rangle (\mod{2p_1})}
\chi^{i}(\gamma)-\frac{p_1-1}{2}p^h \right)\\
=&A(-1)^a\left(-\left(\frac{1-p_1+2\sqrt{-p_1}}{4}\right)\left(\frac{-1+\eta(a)\sqrt{-p_1}}{2}\right)\left(\frac{-1-\sqrt{-p_1}}{2}\right)\right.\\
&\hspace{2cm}-\left.\left(-\frac{1-p_1-2\sqrt{-p_1}}{4}\right)
\left(\frac{-1-\eta(a)\sqrt{-p_1}}{2}\right)\left(\frac{-1+\sqrt{-p_1}}{2}\right)-\frac{p_1-1}{2}p^h \right)\\
=&A(-1)^a\cdot \left\{
\begin{array}{ll}
0&  \mbox{if $a \in \langle p\rangle\,(\mod{p_1})$,}\\
-\frac{p_1^2-6p_1+1}{8}-\frac{p_1-1}{2}p^h&  \mbox{if $a \in -\langle p\rangle\,(\mod{p_1})$.}
 \end{array}
\right.
\end{align*}
This shows that  $\psi(\gamma^a S_2)$
 is  constant according as $\gamma^a\in S_i$, $1\le i\le 5$.
Thus, we have the assertion  of the theorem.
\qed
\vspace{0.3cm}

\begin{remark}\label{examp1}
\begin{itemize}
\item[(i)] 
By using Magma \cite{mamga}, we find that
$(G,\{S_i\}_{i=0}^5)$ becomes a five-class association scheme in the following cases:
\[
(p,p_1,h)=(3,11,1),(5,19,1),(17,67,1),(3,107,3),(41,163,1),(5,499,3).
\]
Note that the association schemes for $(p,p_1)=(3,11),(3,107)$ are nonsymmetric
with exactly two pairs of nonsymmetric relations, where $-S_3=S_2$ and $-S_5=S_4$, and the others are symmetric.
\item[(ii)] We obtain infinite families of five-class association schemes by using the ``recursive'' technique applied in \cite{FX111,FMX11,Momi}: 
Let $p$ and $p_1$ be primes such that $[\Z_{2p_1^m}^\ast:\langle p\rangle]=2$ for all $m\ge 1$, and let
$q=p^{\frac{p_1-1}{2}p_1^{m-1}}$. Define
\begin{align*}
&R_0=\{0\},\, \, R_1=\bigcup_{i=0}^{p_1^{m-1}-1}\bigcup_{j \in \langle p\rangle (\mod{p_1})}C_{2i+p_1^{m-1} j}^{(p_1^m,q)},\\
&\, R_2=\bigcup_{i=0}^{p_1^{m-1}-1}\bigcup_{j \in -\langle p\rangle (\mod{p_1})}C_{2i+p_1^{m-1}j}^{(p_1^m,q)},\, \,
R_3=\bigcup_{i=0}^{p_1^{m-1}-1}C_{2i}^{(p_1^m,q)}, \, \,
\end{align*}
and
\begin{align*}
&S_0=R_0, \, S_1=R_1,\, S_2=\bigcup_{i=0}^{p_1^{m-1}-1}\bigcup_{j \in -\langle p\rangle (\mod{2p_1})}C_{2i+p_1^{m-1}j}^{(2p_1^m,q)}, \, \\
&
S_3=\bigcup_{i=0}^{p_1^{m-1}-1}\bigcup_{j \in -\langle p\rangle (\mod{2p_1})}C_{2i+p_1^{m-1}j+p_1^m}^{(2p_1^m,q)},\, \, 
S_4=\bigcup_{i=0}^{p_1^{m-1}-1}C_{2i}^{(2p_1^m,q)},\, \, 
S_5=\bigcup_{i=0}^{p_1^{m-1}-1}C_{2i+p_1^{m}}^{(2p_1^m,q)}.
\end{align*}
Write $G'=(\F_{q},+)$.
Then, the partition $(G',\{R_i\}_{i=0}^{3})$ becomes  a three-class association scheme and the partition $(G',\{S_i\}_{i=0}^{5})$ forms a  five-class fission scheme of $(G',\{R_i\}_{i=0}^{3})$. 
The proof is done by the quite similar argument as \cite[Theorem 3.6]{FMX11}, that is, reducing the proof  to the case $m=1$, i.e., the proof of Theorem~\ref{thmind22}. 
So we omit the proof. 

Since the pairs
\[
(p,p_1)=(5,19),(17,67),(3,107),(41,163),(5,499)
\]
as listed in Remark~\ref{examp1} (i) satisfy that $[\Z_{2p_1^m}^\ast:\langle p\rangle]=2$ for all $m\ge 1$, we obtain five infinite families of five-class association schemes.  Since $[\Z_{2p_1^m}^\ast:\langle p\rangle]\not=2$ for $m\ge 2$ in the case $(p,p_1)=(3,11)$,  it can not be generalized in this way.
\item[(iii)]
We show that the five infinite families of association schemes of Remark~\ref{examp1} (ii) are
primitive similar to the argument in Remark~\ref{primitivity1}: if $(p,p_1)=(3,107)$, i.e., the schemes are nonsymmetric, by using the  computation of $\psi(\gamma^a R_3)$ given in \cite{FX111} we have
\begin{align*}
\psi(\gamma^a(S_4\cup -S_5))&=\psi(\gamma^a R_3)\\&=
\frac{1}{p_1}\left(p^{(f-h)/2}\left(\frac{-b\pm p_1c}{2}\right)-1\right)\, \, \mbox{ or }\,\, \frac{1}{p_1}\left(p^{(f-h)/2}\left(\frac{p_1-1}{2}\right)-1\right)\\
&=\frac{3^{\frac{f-1}{2}}\cdot 53-1}{107} \, \mbox{ or } \frac{-3^{\frac{f-1}{2}}\cdot 54-1}{107} \not=\frac{3^f-1}{107}.
\end{align*}
Next consider the case $(p,p_1)\not=(3,107)$, i.e., the schemes are symmetric. By \cite[Theorem 3.2]{Momi}, if $p\equiv 1\,(\mod{4})$, it holds that for any $I\subseteq \Z_{2p_1}$ \begin{equation}\label{recur}
\sum_{i=0}^{p_1^m-1}\sum_{j\in I}\psi(\gamma^a C_{2i+p_1^{m} j}^{(2p_1^{m+1},p^{\frac{p_1-1}{2}p_1^m})})= p^{\frac{(p_1-1)^2p_1^{m-1}}{4}}
\sum_{i=0}^{p_1^{m-1}-1}\sum_{j\in I}\psi'(\omega^a C_{2i+p_1^{m-1}j}^{(2p_1^{m},p^{\frac{p_1-1}{2}p_1^{m-1}})})+\frac{p^{\frac{(p_1-1)^2p_1^{m-1}}{4}}-1}{p_1},
\end{equation}
where  $\psi$ and $\psi'$ are
canonical additive characters of $\F_{p^{\frac{p_1-1}{2}p_1^{m}}}$ and $\F_{p^{\frac{p_1-1}{2}p_1^{m-1}}}$, and  $\gamma$ and $\omega$ are
primitive roots of $\F_{p^{\frac{p_1-1}{2}p_1^{m}}}$ and $\F_{p^{\frac{p_1-1}{2}p_1^{m-1}}}$, respectively. Then, by the computation of $\psi(\gamma^a S_4)$ for $m=1$ in Theorem~\ref{thmind22}, we have that for general $m$ 
\begin{align*}
\psi(\gamma^a S_4)=&\frac{1}{2p_1}\left(p^{(f-h)/2}\left(\frac{-b-p_1c}{2}\right)-1\right),\\
&\, \, \frac{1}{2p_1}\left(p^{(f-h)/2}\left(\frac{-b+p_1c}{2}\right)-1\right)\pm \frac{p^{\frac{f}{2}-h}}{2p_1}\left(p^h+\frac{3p_1-1}{4}\right),\\
&\, \frac{1}{2p_1}\left(p^{(f-h)/2}\left(\frac{p_1-1}{2}\right)-1\right)\pm \frac{p^{\frac{f}{2}-h}}{2p_1}\left(p^h-\frac{(p_1-1)^2}{4}\right),
\end{align*}
which are not equal to  $\frac{p^f-1}{2p_1}$ for any $(p,p_1)=(5,19),(17,67),(41,163),(5,499)$.
\end{itemize}
\end{remark}
\section{A nonsymmetric four-class fission scheme  of a conference graph}
First we recall the construction of two-class translation schemes in \cite{fx}.  A more common way to describe two-class translation schemes is to use the language of partial difference sets, but we don't introduce this definition here and refer to \cite{pds_sur} for details.  Let $p_1$ be an odd prime. We assume the following specific index 2 case: $N_0=2p_1 $, $p_1>3$ is a prime, and $p_1\equiv 3\pmod{4}$; $p$ is an odd prime such that $[\Z_{N_0}^*:\langle p\rangle]=2$ (that is, $f:={\rm ord}_{N_0}(p)=\phi(N_0)/2$). Assume that $p_1\equiv 7\pmod{8}$. Let $C_i=C_{i}^{(N_0,q)}$ for $1\leq i\leq N_0-1$.

\begin{thm}\label{7mod8}
Assume that we are in the index 2 case specified as above. Let $p_1\equiv 7\pmod{8}$, $f=(p_1-1)/2$, and $q=p^f$, and write $G=(\F_q,+)$. Let $I_0$ be any subset of $\Z_{N_0}$ such that $\{i\pmod{p_1}\mid i \in I_0\}=\Z_{p_1}$,  and let $D_0=\bigcup_{i\in I_0}C_i^{(N_0,q)}$, $D_1=\F_q^*\setminus D_0$. Then $(G,\{\{0\},D_0,D_1\})$ is a two-class translation scheme if $p\equiv 1\pmod{4}$.
\end{thm}

The Cayley graph of $D_0$ in Theorem \ref{7mod8} is a conference graph, namely, it is a strongly regular graph with parameters of the form $(v,\frac{v-1}{2},\frac{v-5}{4},\frac{v-1}{4})$.
By using Magma \cite{mamga}, we find that when $p=37$, $q=p^3$, and $N=4p_1=28$, the index sets
\[
I_1=\{ 0, 1, 4, 12, 16, 20, 24 \}\subset\Z_{28},\, I_2=I_1+7,\, I_3=I_1+14,\,I_4=I_1+21
\]
define a four-class nonsymmetric scheme whose symmetrization is covered by Theorem \ref{7mod8}. The dual scheme is defined by
\begin{align*}
&J_1=\{ 0, 4, 8, 12, 13, 16, 24 \}\subset\Z_{28},\,
J_2=J_1+7,\, J_3=J_1+14, \,J_4=J_1+21.
\end{align*}

We have $q=37^3=107^2+4\cdot99^2=37^2+4\cdot111^2$. By Ma and Wang's results \cite{mw}, the putative skew-symmetric fission scheme with $4$ classes of a conference graph exists only if $q\equiv 5\pmod{8}$ and there exist integers $g,h$ such that $q=g^2+4h^2$ with $g\equiv 1\pmod{4}$. Moreover, the first eigenmatrix $P$ is of the form
\[
P=\begin{pmatrix}
1 &f & f & f &f\\
1 &\rho & \tau & \bar{\rho} & \bar{\tau} \\
1 &\tau & \bar{\rho} & \bar{\tau} & \rho \\
1 & \bar{\rho} & \bar{\tau}  & \rho  & \tau \\
1 & \bar{\tau} & \rho & \tau &  \bar{\rho}
\end{pmatrix}
\]
where $f=\frac{q-1}{4}$, $\rho=\frac{1}{4}(-1+\sqrt{q}+\sqrt{-2q+2g\sqrt{q}})$,
$\tau=\frac{1}{4}(-1-\sqrt{q}+\sqrt{-2q-2g\sqrt{q}})$. It turns out that our example happens with $g=37$.  For the corresponding cyclotomic scheme of order $4$ over $\F_{37^3}$, we have $g=-107$.

The intersection matrices of our scheme are given below:
\begin{align*}
 B_1=\begin{pmatrix}    0&     0   &  0 &12663  &   0 \\
     1 & 3170 & 3161 & 3170 & 3161 \\
     0  &3115&  3161 & 3161&  3226 \\
     0  &3152 & 3226 & 3170 & 3115 \\
     0  &3226 & 3115 & 3161 & 3161
    \end{pmatrix}&&
B_2=\begin{pmatrix}
     0  &   0   &  0 &    0& 12663 \\
     0 & 3161 & 3226&  3115  &3161 \\
     1 & 3161  &3170 & 3161  &3170 \\
     0 & 3226 & 3115&  3161 & 3161 \\
     0  &3115 & 3152 & 3226 & 3170
     \end{pmatrix}\\
B_3=\begin{pmatrix}
     0 &12663  &   0  &   0  &   0 \\
     0 & 3170 & 3115 & 3152&  3226 \\
     0 & 3161 & 3161 & 3226 & 3115 \\
     1 & 3170 & 3161 & 3170 & 3161 \\
     0 & 3161 & 3226&  3115 & 3161
    \end{pmatrix}&&
B_4=\begin{pmatrix}
     0   &  0 &12663  &   0  &   0 \\
     0  &3161&  3161  &3226&  3115 \\
     0 & 3226  &3170 & 3115 & 3152 \\
     0  &3115 & 3161 & 3161 & 3226 \\
     1 & 3161 & 3170 & 3161 & 3170
     \end{pmatrix}
\end{align*}

\begin{remark} The cyclotomic field $\Q(\xi_{37})$ has a cyclic Galois group so has only one subfield of extension degree $4$ over $\Q$, whose algebraic integer ring has  a basis given by the Gauss periods of index $4$ over $\F_{37}$, denoted by $\eta_0,\ldots,\eta_3$ as usual. We can express the $\rho$, $\tau$   in term of the linear combinations of these $\eta_i$'s over $\Z$ as follows: with $g=37$, we have $\rho=9+37\eta_0$; with $g=-107$, we have $\rho=8+35\eta_0+5\eta_1-7\eta_3$.
\end{remark}

\end{document}